\documentstyle[amssymb]{amsart}
 \newcommand{\inputpath}{}

\newcommand{\stopproof}{\hfill \nobreak\medskip $\blacksquare$ \\
\hspace*{\fill}}

\newcommand{\dom}{\mbox{\rm dom}}

\newcommand{\AND}{\mbox{ \rm and }}

\newcommand{\IF}{\mbox{ \rm if }}

\newcommand{\OR}{\mbox{ \rm or }}

\newcommand{\forces}[2]{\Vdash_{#1} \mbox{``} #2 \mbox{''}}

\newcommand{\proof}{{\bf Proof:} \ }

\newcommand{\PP}{{\Bbb P}}

\newcommand{\RR}{{\Bbb R}}

\newcommand{\KK}{{\Bbb K}}

\newcommand{\SS}{{\Bbb S}}

\newcommand{\presup}[2]{\, ^{#1} \! #2}

\newcommand{\lcard}{\, \mid \!}
\newcommand{\rcard}{\! \mid \,}

\newcommand{\fomom}{\presup{\omega}{\omega}}
\newcommand{\wfomom}{\presup{\stackrel{\omega}{\smile}}{\omega}}

\newcommand{\card}[1]{\lcard #1 \rcard}

\newtheorem{defin}{Definition}[section]

\newtheorem{lemma}{Lemma}[section]

\newcommand{\Cantor}{{[0,1]}}

\newcommand{\cov}[1]{\mbox{cov}(#1)}
\newcommand{\RestrictedTree}[1]{\langle #1\rangle}
\newcommand{\SuperPerfect}{\SS}
\newcommand{\wpresup}[2]{\presup{\stackrel{#1}{\smile}}{#2}}
\title{Decomposing Baire Class 1 Functions into Continuous Functions}
\author{Saharon Shelah}
\address{Institute of Mathematics \\ Hebrew University \\ Jerusalem,
Givat Ram,  Israel and Department of Mathematics\\
Rutgers University\\ New Brunswick, New Jersey}
\author{Juris Stepr\={a}ns}

\address{Department of Mathematics, York University \\
4700 Keele Street \\ North York, Ontario \\ Canada \ \ \ \ M3J 1P3}            \thanks{The first author is partially  supported by the basic reasearch
fund of the Israeli Academy. The second author is partially
supported by NSERC and was a guest of Rutgers University while the
research on this paper was being done.
This is number ??? on the first author's list of publications}
\begin{document}
  \bibliographystyle{\inputpath amsplain}
\maketitle
\section{Introduction}
In \cite{cimopaso} the authors consider the following question:
What is the least cardinal $\kappa$ such that every function of first
Baire class can be decomposed into $\kappa$ continuous functions? This
cardinal $\kappa$ will be denoted by $\frak{dec}$. The authors of
\cite{cimopaso} were able to show that $\cov{\KK} \leq \frak{dec}
\leq \frak d$ and asked whether these inequalities could,
consistently, be strict. By $\cov{\KK} $ is meant the least number of
closed nowhere dense sets required to cover the real line and by
$\frak d$ is denoted the least cardinal of a dominating family in $\fomom$.
In \cite{step.28} it was shown that it is consistent that
 $\cov{\KK} \neq \frak{dec}$. In this paper it will be shown that the
second inequality can also be made strict. The model where ${\frak dec}$
is different from $\frak d$ is the one obtained by adding $\omega_2$
Miller --- sometimes known as super-perfect or rational-perfect --- reals to a model of the
Continuum Hypothesis. It is somewhat surprising that the model used
to establish the consistency of the other inequality, $\cov{\KK} \neq \frak{dec}$, is a
slight modification of the iteration of super-perfect forcing.

By ${\presup{\stackrel{\omega}{\smile}}{\omega}} $ will be denoted
$\cup_{n\in\omega}\{\presup{n}{\omega} : n\in\omega\}$.
As usual, a tree will be defined to mean an initial subset of
${\presup{\stackrel{\omega}{\smile}}{\omega}} $ under
$\subseteq$.
 So if $T$ is a tree and $t\in T$ then $t\restriction k \in T$ for
each $k\in\omega$.
Also,  $T\RestrictedTree{t}$ will be defined to be
$\{s \in T: s\subseteq t \OR
t\subseteq s\}$.
If  $t$  and $s$ are both finite sequences
 then $s\wedge t$ is defined by declaring that $\dom(s \wedge t ) =
\card{\dom(t)} + \card{\dom(s)}$ and
$$s\wedge t (i) = \left\{\begin{array}{ll}
                     s(i)  & \IF i\in\dom(s)  \\
                      t(i \  - \  \card{\dom(s)})  & \IF i\notin\dom(s)
                  \end{array}\right.$$
If $t\in T \subseteq  \wfomom $
then $i\in \omega$ then $t\wedge i$ is defined to be
$t \wedge \{(0,i)\} $
and           $i\wedge t $ is defined to be
$\{(0,i)\}\wedge t $.
                  Finally,
 $\overline{T} = \{f\in
\fomom
 : (\forall n\in\omega)
(f\restriction n \in T)\}$ and closure in other spaces is denoted similarly.

\begin{defin}
If $T\subseteq \wfomom$ is a tree then $\beta(T)$ will be defined to
be the set of all $t\in T$ such that $\card{\{n \in\omega : t\wedge
n\in T\}} = \aleph_0$. A tree $T\subseteq \wfomom$ is said to be {\em
super-perfect} if for each $t\in T$ there is some $s\in \beta(T)$ such
that $t\subseteq s$ and if
 $\card{\{n \in\omega : t\wedge
n\in T\}} \in\{1,\aleph_0\}$ for each $t\in T$. The set of all
super-perfect
 trees will be denoted
by $\SuperPerfect$.
\end{defin}

For each $T\in \SuperPerfect$ there is a natural way to assign a mapping
$\theta_T :\wfomom\to\beta(T)$  such that:
\begin{itemize}
\item $\theta_T$ is one-to-one and onto $\beta(T)$
\item  $s\subseteq t$ if and only if $\theta(s)\subseteq\theta(t)$
\item   $s\leq_{\mbox{Lex}} t$ if and only if $\theta(s)
\leq_{\mbox{Lex}}
\theta(t)$.
\end{itemize}
Notice that $\theta_T(\emptyset)$ is the root of $T$.
Using the mapping $\theta_T$, it is possible to define a  refinement of the
ordering on $\SuperPerfect$.
\begin{defin}
Define $T\prec_n S$ if both $S$ and $T$ are in $\SuperPerfect$, $T\subseteq S$
and $\theta_T\restriction \presup{n}{\omega} =
\theta_S\restriction \presup{n}{\omega}$.
\end{defin}
It should be clear that the ordering $\prec_n$ satisfies Axiom A.
The proof of the main result of this paper will use a fusion
based on a sequence of the orderings $\prec_n$. Notice that while
 $\prec_n$ can be used in the same way as the analogous ordering for
Sacks reals in the case of adding a single real,
is not as easy to deal with in the context of iterations.
The chief difficulty is that
 $\prec_n$ requires deciding an infinite amount of information because
branching is infinite. This conflicts with the usual goal of fusion
arguments which decide only a finite amount of information at a time.

 \section{Iterated Super-Perfect Reals}
It will be shown that in a model obtained by iterating $\omega_2$ times the partial
orders $\SuperPerfect$ with countable
support over a ground model where
$2^{\aleph_0} = \aleph_1$ yields a model where $\frak{d} = \aleph_2$
and $\frak{dec} = \aleph_1$. The fact that  $\frak{d} = \aleph_2$ is
well known \cite{mill.ratpfct}. The fact that  $\frak{dec} = \aleph_1$
is an immediate consequence of the following result.
\begin{lemma}
Suppose that $\xi\in\omega_2+1$, $\SuperPerfect_\xi$ is the  iteration with countable
support of the partial orders $\SuperPerfect$ and\label{keylem} $G$  is
$\SuperPerfect_\xi$-generic over $V$. Then for any $x\in \Cantor$ in $V[G]$
 and any Borel function
$H:\Cantor\to\Cantor$ in $V[G]$ there is a Borel set $X\in V$ such that $x\in
X$ and $H\restriction X$ is continuous.
\end{lemma}

Saying that $X \in V$ means, of course, that the real coding the
Borel set $X$ belongs to the model $V$.
In order to prove Lemma~\ref{keylem} it will be useful to employ a different
interpretation of iterated super-perfect forcing.  The next sequence of
definitions will be used in doing this.  If $G$ is $\SuperPerfect_\xi$-generic
over some model $\frak M$ then there is a natural way to assign a mapping
$\Gamma :  \xi \cap {\frak M} \to \fomom$ such that ${\frak M}[G] = {\frak
M}[\Gamma]$.  On the other hand, given $\Gamma :{\frak M}\cap \xi \to
\fomom$ define
$G_\Gamma({\frak M})$ to be the set
 $$\{q\in {\frak M}\cap \SuperPerfect_\xi : \forall k \in
\omega\forall A\in [{\frak M}\cap\xi]^{<\aleph_0}\exists p \leq q
\forall \alpha \in A(p\restriction \alpha
\forces{\SuperPerfect_\alpha}{\Gamma(\alpha) \restriction k\in
p(\alpha)})\}$$
and say that $\Gamma$ is $\SuperPerfect_\xi$-generic over ${\frak M}$
if and only if $G_\Gamma$ is $\SuperPerfect_\xi$-generic over ${\frak M}$.
Note that if $G$ is $\SuperPerfect_\xi$-generic over ${\frak M}$ and
$\Gamma:
{\frak M}\cap \xi \to \fomom$ is its associated function then
$G_\Gamma({\frak M}) = G$.
This will be used without further comment to identify
$\SuperPerfect_\xi$-generic sets over $\frak M$ with elements of
$(\fomom)^{{\frak M}\cap \xi}$.  Whenever a
topology on $(\fomom)^X$ is mentioned, the product topology is intended.

\begin{defin}
If $p \in \SuperPerfect_\xi$ and $\Lambda \in [\xi]^{\leq\aleph_0}$
then define $S(\Lambda,p) $ to be the set of all functions $\Gamma :
\Lambda \to \fomom$ such that for all $k\in \omega$ and for all finite
subsets $A\subseteq \Lambda$ there is $q \leq p$ such
that $q\forces{\SuperPerfect_\xi}{\Gamma(\alpha)\restriction k \in
q(\alpha)}$
for all $\alpha\in A$.
 \end{defin}

\begin{defin} Given a countable elementary submodel ${\frak M} \prec H((2^{\aleph_0})^+)$
and $p \in \SuperPerfect_\xi$ define $p$ to be strongly\label{StronglyGeneric}
$\SuperPerfect_\xi$-generic over $\frak M$ if and only if \begin{itemize} \item
each $\Gamma \in S({\frak M}\cap \xi,p)$ is $\SuperPerfect_\xi$-generic over
$\frak M$ \item if $\psi$ is a statement of the $\SuperPerfect_\xi$-forcing
language using only parameters from $\frak M$, then $\{\Gamma \in S({\frak
M}\cap \xi,p)\ : \ {\frak M}[\Gamma] \models \psi\}$ is a clopen set
in $S({\frak M}\cap \xi,p)$.
\end{itemize} \end{defin}
A set $X \subseteq (\fomom)^\alpha$ will be defined to be {\em large} by
induction on $\alpha$.
\begin{defin}
  If $\alpha = 1$ then $X$ is large if $X$ is a
superperfect tree.  If $\alpha$ is a limit then $X$ is large if the projection
of $X$ to $(\fomom)^\beta$ is large for every $\beta\in\alpha$.  If $\alpha =
\beta + 1$ then $X$ is large if there is a large set $Y\subseteq (\fomom)^\beta
$ such that $X = \cup_{y\in Y}\{y\}\times X_y$ and each $X_y$ is a large subset
of $\fomom$.  \end{defin}

>From large closed sets it is possible to obtain, in a
natural way, conditions in $\SuperPerfect_\xi$.  \begin{defin} If $X \subseteq
(\fomom)^\alpha$ is a large closed set then define $p_X\in \SuperPerfect_\alpha$
by defining $p_X(\eta)$ to be the $\SuperPerfect_\eta $ name for that subset
$T\subseteq \fomom$ such that if $\Gamma:\alpha \to \fomom$ is $\SuperPerfect_\alpha$
generic then
$$T = \{ f\in \fomom | (\exists h)(\Gamma\restriction \eta \cup \{(\eta,f)\}
\cup h\in X)\}$$\end{defin} Observe that, if $X\subseteq (\fomom)^\alpha$ is
large and closed, it follows that $p_X\in \SuperPerfect_\alpha$. The
following result provides a partial converse to this observation.

\begin{lemma}
If $p \in \SuperPerfect_\xi$ and ${\frak M} \prec H((2^{\aleph_0})^+)$ is a
countable elementary submodel containing $p$ then there is \label{ShelCountSup}
$q \leq p$ such that $q$ is strongly $\SuperPerfect_\xi$-generic over
$\frak M$.
\end{lemma}
\proof The proof consists of merely
repeating the proof that the countable support iteration of proper
partial orders is
proper and checking the assertions in this special case. Only a sketch will be
given and the reader should consult \cite{shel.pf} for details.

The proof  is by induction on $\xi$. If $\xi = 1$ then a standard fusion
argument
applied to an enumeration $\{D_n : n\in \omega\}$, of all dense subsets of
$\SuperPerfect$ provides the
result.       In particular, there is a sequence $\{T_i : i\in \omega\}$ such
that $T_{i+1} \prec_i T_i$, $T_0 = T$ and such that           $T_i\langle
\theta_{T_i}(\sigma)\rangle \in D_{i-1}$ for each $\sigma : i \to \omega$. The
condition $T_\omega = \cap_{i\in \omega}T_i$ has the desired property. The fact
that if $\psi$ is a statement of the $\SuperPerfect_\xi$-forcing
language using only parameters from $\frak M$, then $\{\Gamma \in S({\frak
M},T_\omega)\ : \ {\frak M}[\Gamma] \models \psi\}$ is a clopen
set is obvious because $S(1, T_\omega) = \overline{T_\omega}$.

If $\xi = \mu + 1$ then use the induction hypothesis to find $q' \leq
p\restriction \xi$ such that $q'$ is strongly $\SuperPerfect_\mu$-generic over
$\frak M$. Then, in particular, $q'$ is  $\SuperPerfect_\mu$-generic over
$\frak M$     and so, if $G$ contains $q'$ and is  $\SuperPerfect_\mu$-generic
over $V$ it is also generic over $\frak M$. Therefore ${\frak M}[G] $ is an
elementary submodel in $V[G]$ and it is possible to choose an enumeration $\{D_n
: n\in\omega\}$ of all dense subsets of $\SuperPerfect$ which are members
of ${\frak M}[G] $. It is therefore possible to choose, in ${\frak M}[G]$, as in the case $\xi =
1$, a sequence  $\{T_i : i\in \omega\}$ such
that $T_{i+1} \prec_i T_i$ and such that           $T_i\langle
\theta_{T_i}(\sigma)\rangle \in D_{i-1}$ for each $\sigma : i \to \omega$. The
condition $T_{\omega} = \cap_{i\in \omega}T_i$ is then strongly $\SuperPerfect$-generic over ${\frak
M}[G]$. Notice that, while $T_{\omega}$
does not itself have a name in $\frak M$, each $T_n$ does have a name
and so there are enough objects in ${\frak M}[G]$ to construct $T_{\omega}$.

In order to see that $q=q'  * T_{\omega}$
is strongly $\SuperPerfect_\xi$-generic over ${\frak
M}$ suppose that $\Gamma \in S({\frak M}\cap
\xi,q)$. Obvioulsy $\Gamma\restriction \mu \in S({\frak M}\cap
\mu,q')$                  and therefore ${\frak M}[\Gamma]$ is an elementary
submodel.  Hence, by
genericity,  $T_{i+1} \prec_i T_i$, $T_0 = T$ and $T_i\langle
\theta_{T_i}(\sigma)\rangle \in D_{i-1}$ and so it follows that
$\cap\{T_i : i\in \omega\}$ is a strongly $\SuperPerfect$-generic
condition over ${\frak M}[G]$. Hence $\Gamma(\xi)$ is  $\SuperPerfect$-generic
over  ${\frak M}[G]$ and so $\Gamma$ is  $\SuperPerfect_\xi$-generic
over $\frak M$.

Just as in the case $\xi = 1$, it is easy to use the induction hypothesis to
see that
 if $\psi$ is a statement of the $\SuperPerfect_\xi$-forcing
language using only parameters from $\frak M$, then $\{\Gamma \in S({\frak
M}\cap \xi,q)\ : \ {\frak M}[\Gamma] \models \psi\}$ is a clopen set.

Finally, suppose that $\xi$ is a limit ordinal. If it has uncountable
cofinality then there is nothing to do because of the countable
support of the iteration. So assume that $\{\mu_n : n\in\omega\}$ is an increasing
sequence of ordinals cofinal in $\xi$. Let $\{D_n : n\in \omega\}$
enumerate all dense subsets of $\frak M$ and choose a sequence of
conditions $\{p_i : i\in \omega\}$ such that
\begin{itemize}
\item $p_i\restriction \mu_i$ is strongly $\SuperPerfect_{\mu_i}$-generic over
$\frak M$
\item  $p\restriction \mu_i\forces{\SuperPerfect_{\mu_i}}{p_i\restriction (\xi
\setminus \mu_i) \in D_i/G}$ (this is an abbreviation for the  more
precise statement:
$$p\restriction \mu_i\forces{\SuperPerfect_{\mu_i}}{(\exists q\in G\cap
\SuperPerfect_{\mu_i})(
q * p_i\restriction (\xi
\setminus \mu_i) \in D_i)}$$ and will be used later as well)
\item $p_i\restriction (\xi \setminus \mu_i)$ belongs to $\frak M$
\item $p\restriction \mu_i\forces{\SuperPerfect_{\mu_i}}{p_{i+1}\restriction
(\mu_{i+1} \setminus \mu_i) \mbox{ is }
\SuperPerfect_{\mu_{i+1}\setminus \mu_i}\mbox{-generic over }{\frak
M}[G]}$
\item $p_{i+1} \leq p_i$
\end{itemize}
Notice that the statement that ${p_i\restriction (\xi
\setminus \mu_i) \in D_i/G}$ can be expressed in $\frak M$ and so if $\Gamma
\in S({\frak M}\cap \SuperPerfect_{\mu_i}, p_i\restriction \mu_i)$ then
$p_i\restriction (\xi
\setminus \mu_i) \in D_i/\Gamma$. From this it easily follows that letting
$p_\omega = \lim_{n\in\omega}p_n$ yields a strongly $\SuperPerfect_\xi$-generic
condition over $\frak M$.

To see that
 if $\psi$ is a statement of the $\SuperPerfect_\xi$-forcing
language using only parameters from $\frak M$, then $\{\Gamma \in S({\frak
M}\cap \xi,p_\omega)\ : \ {\frak M}[\Gamma] \models \psi\}$ is a clopen set,
observe that for any such
$\psi$ there corresponds the dense subset of $\SuperPerfect_\xi$ consisting of
all
conditions which decide $\psi$. Any such dense set is therefore $D_n$ for some
$n\in \omega$.    It follows that if $\Gamma \in S({\frak
M}\cap \xi,p_\omega)$ then the interpretation of $p_n\restriction(\xi
\setminus\mu_n)$ in ${\frak M}[\Gamma\restriction \mu_n]$ decides the truth
value of $\psi$ because $p_n\restriction \mu_n$ is strongly
$\SuperPerfect_{\mu_n}$-generic over ${\frak M}$.
>From the induction hypothesis it follows that there is a clopen set
$U \subseteq S({\frak M}\cap \mu_n, p_n\restriction \mu_n)$ such that for each
$\Gamma' \in U$ the model ${\frak M}[\Gamma']$ satisfies that
the interpretation of $p_n\restriction(\xi
\setminus\mu_n)$ in ${\frak M}[\Gamma\restriction \mu_n]$
 decides the truth value of
$\psi$. Let $U^*$ be the lifting of $U$ to
$S({\frak M}\cap \xi, p_{\omega})$ --- in other words, $\Gamma \in
U^*$ if and only if $\Gamma\restriction \mu_n \in U$. Since
the interpretation of $p_\omega\restriction(\xi
\setminus\mu_n)$ in ${\frak M}[\Gamma\restriction \mu_n]$ is a
stronger condition than the interpretation of $p_n\restriction(\xi
\setminus\mu_n)$ in ${\frak M}[\Gamma\restriction \mu_n]$, it follows
that  $U^* \subseteq S({\frak
M}\cap \xi,p_\omega)$ is the desired clopen set.
\stopproof

\begin{defin}
           A subset $X\subseteq \presup{n}{\omega}$ is said to be a {\em full subset} if, $X\neq \emptyset$
           and for each $x\in X$ and $i \in n$ there is $A\in [\omega]^{\aleph_0}$ such that for all $m \in A$ there is
           $x_m \in X$ such that $x_m \restriction i = x\restriction i$ and $x_m(i) = m$.
\end{defin}
\begin{lemma}
                 If $F: \presup{n}{\omega} \to [0,1]$ is \label{FullSubtree} a one-to-one function
                 then there is a full subset $T\subseteq   \presup{n}{\omega}$ such that the image of
                 $T$ under $F$ is discrete.
      \end{lemma}
\proof Proceed by induction on $n$ to prove the following stronger assertion: If $F: \presup{n}{\omega}
\to [0,1]$ is one-to-one then there is a full subset $T\subseteq \presup{n}{\omega}$, there is $f\in \fomom$
and there is $x\in [0,1]$
such that \begin{itemize}
\item[{\bf A.}] for any  a descending sequence $\{U_i : i\in \omega\}$
 of neighbourhoods of $x$ such that $\mbox{diam}(U_{n+1})\cdot f(\lceil 1/\mbox{diam}
(U_n)\rceil ) < 1$                            and for each $X\in [\omega]^{\aleph_0}$
the set $\{t\in T : F(t) \in \cup_{i\in X } (U_i\setminus \overline{U_{i+1}})\}$
is a full subset.
           \end{itemize}

The case $n=1$ is easy. Choose $A\in [\omega]^{\aleph_0}$ such that $\{F(\emptyset \wedge i) : i \in A\}$
converges to $x\in [0,1]$. Let $f \in \fomom$ be any increasing function such that for each  $m\in\omega$
there is some $j\in A$ such that $1/m > | F(\emptyset \wedge j) | >
1/f(m)$.
 Let $T = \{\emptyset \wedge i : i\in A\}$.

Now let $F:\presup{n+1}{\omega} \to [0,1]$ be one-to-one.  Use the induction
hypothesis to find, for each $m\in\omega$, full subsets $T_m\subseteq
\presup{n}{\omega}$ such that the
image of $F$ restricted to $$\{x\in \presup{n+1}{\omega} : (\exists t\in T_m)(x=
\emptyset\wedge  m\wedge t)\}$$ is a discrete family and Condition {\bf A.}
is witnessed by $f_m\in \fomom$ and $x_m \in [0,1]$.  There are two cases to
consider depending on whether or not there is $Z\in [\omega]^{\aleph_0}$ such
that $\{x_m | m\in Z\}$ are all distinct.

\noindent{\bf Case 1}

Assume that there is $Z\in [\omega]^{\aleph_0}$ such that $\{x_m : m\in Z\}$ are all distinct. It is
then possible to assume that there is some $x\in [0,1]$ such that $\lim_{n\in Z}x_m = x$ and that,
without loss of generality,
$x_m > x_{m+1} > x$. As in the case $n=1$, it is possible to find $f\in\fomom$ such that
for any  a descending sequence  $\{U_i : i\in \omega\}$
 of neighbourhoods of $x$ such that $\mbox{diam}(U_{n+1})\cdot f(\lceil 1/\mbox{diam}
(U_n)\rceil ) < 1$ and for each  $X \in [\omega]^{\aleph_0}$ the set $\{
m\in\omega : x_m \in \cup_{i\in X } (U_i\setminus \overline{
U_{i+1}})\}$
is infinite. Notice that each $U_i\setminus \overline{
U_{i+1}}$ is open, so it follows from Condition {\bf A.} that $\{t\in T_m
: F( m \wedge t) \in U_i\setminus \overline{
U_{i+1}}\}$ is a full subset provided that $ x_m \in U_i\setminus \overline{
U_{i+1}}$. Hence, $$\cup\{ \{t\in T_m
: F(\langle m \rangle\wedge t) \in U_i\setminus \overline{
U_{i+1}}\}  : x_m \in U_i\setminus \overline{
U_{i+1}}      \}$$ is a full subset provided that
$\mbox{diam}(U_{n+1})\cdot f(\lceil 1/\mbox{diam}
(U_n)\rceil ) < 1$ and  $X \in [\omega]^{\aleph_0}$.
Let $T = \{t\in \presup{n+1}{\omega} : (\exists t'\in
T_{t(0)})(t= t(0)\wedge t')\}$. Then $T$,
$f$ and $x$ satisfy the Condition {\bf A.}

\noindent{\bf Case 2}

In this case there exists $x\in [0,1]$ such
that $x_m = x$ for all but finitely many $m \in \omega$. Let $f\in \fomom$ be such that $f \geq^* f_m$ for all $m \in \omega$.
Let $$T = \{t\in \presup{n+1}{\omega} : (\exists t'\in T_{t(0)})(t=
t(0)\wedge t' \AND
x_{t(0)} = x)\}$$ To see that this works, suppose that $\{U_i : i\in \omega\}$
is a descending sequence of neighbourhoods of $x$
such that $\mbox{diam}(U_{i+1})\cdot f(\lceil 1/\mbox{diam}
(U_i)\rceil) < 1$ and suppose that $X\in [\omega]^{\aleph_0}$.

Let $X = \cup_{j\in\omega}X_j$ be a partition of $X$ into infinite
subsets. It may be assumed that $f(i) \geq f_m(i)$ for all $i \in
X_m$.
 By the induction hypothesis it follows that $\{t\in T_m :
F(t) \in \cup_{i\in X_m } (U_i\setminus \overline{U_{i+1}})\}$ is a full subset
of $\presup{n}{\omega}$ for each $m \in \omega$ because $f \geq^*
f_m$.  Hence $\{t\in T : F(t) \in \cup_{i\in X } (U_i\setminus
\overline{U_{i+1}})\}$ is a full subset of $\presup{n+1}{\omega}$.  \stopproof

Although this fact will not be used, it should be noted that
Lemma~\ref{FullSubtree} can be generalised to arbitrary well founded trees.

If $X\subseteq (\fomom)^\alpha$ is large then for each $e :
\beta \to \fomom$ let $X_e$ represent the set of all $f :
{\alpha \setminus \beta}\to \fomom$ such that $e\cup f \in X$.
Note that if $ h \in X$ then for every $\beta \in \alpha$,
$X_{h\restriction \beta}$ is a large subset of
$(\fomom)^{\alpha\setminus \beta}$.  Moreover, the projection
$X_{h\restriction \beta}$ to $(\fomom)^{\delta\setminus \beta}$ is
large provided that $\beta \in \delta$. This set will be denoted by
$\pi_\delta(X_{f\restriction \beta})$.  Note that $\pi_{\beta +
1}(X_{f\restriction \beta})$ is the closure of a super-perfect tree, $T_{X,f,\beta}$ and so
$\theta_{T_{X,f,\beta} } : \wfomom \to T_{X,f,\beta} $ is an isomorphism. This
induces a natural isomorphism from $\wpresup{\alpha}{(\wfomom)} $ to the open sets of $X$ which
will be denoted by $\Phi_X$.

\begin{lemma}
If $\alpha\in \omega_1$, $\frak M$ is a countable elementary submodel, $q \in \SS_\alpha$ and
$F:S({\frak M}\cap\alpha,q)
\to \RR$ is continuous satisfying\begin{itemize} \item[{\bf B.}] for each $\beta\in \alpha$ and each
$e\in (\fomom)^\beta$, if $S({\frak M}\cap\alpha,q)_e\neq \emptyset$, then the range of  $F$
restricted to $S({\frak M}\cap\alpha,q)_e$ is uncountable \end{itemize}
then there is a large                                          \label{GenericCondition}
closed set $X\subseteq S({\frak M}\cap\alpha,q)$ such that $F\restriction X$ is one-to-one and, moreover,
$F\restriction X$ is a homeomorphism onto its range.
\end{lemma}
\proof  For
$\tau \in \wpresup{\alpha}{(\wfomom)}$     and $\tau' \in
\wpresup{\alpha}{(\wfomom)}$                        define $\tau \leq \tau'$ if
and only if $\tau(\sigma) \subseteq \tau'(\sigma)$ for each $\sigma$ in the
domain of $\tau$ and, define $\tau_1$ and $\tau_2$ to be incompatible if there is no $\tau'$
such that $\tau_1 \leq \tau'$ and $\tau_2\leq \tau'$.
To begin, let $\{\tau_i :  i\in\omega\}$ enumerate a subset of $\wpresup{\alpha}{(\wfomom)}$
which forms a tree base for   $S({\frak M}\cap\alpha,q)$ --- in other words, if $i$ and $j$ are in $\omega$
then either
$\tau_i < \tau_j$, $\tau_j < \tau_i$ or $\tau_i$ and $\tau_j$ are incompatible and, moreover,
$\{\Phi_{S({\frak M}\cap\alpha,q)}(\tau_i) :\i\in \omega\}$ is a base for $S({\frak
M}\cap\alpha,q)$. It may also be assumed that
if $\tau_i < \tau_j$ then $i \leq j$ and that
for each $k\in \omega$ there is a unique $\rho$ and some $i\in k$ such that
$\tau_k(\mu) = \tau_i(\mu)$ if $\mu\neq \rho$
and $\tau_k(\rho) = \tau_i(\rho)\wedge W$ for some integer $W$.  Let $X_0 =
S({\frak M}\cap\alpha,q)$. Construct by
induction a sequence $\{(X_k, \{U_i :  i \in k\} :  k\in \omega\}$ such that:
\begin{itemize} \item[{\bf a.}] $X_k$ is a large and closed subset of
$(\fomom)^\alpha$
\item[{\bf b.}]
each $U_i$ is an open subset of $\RR$  \item[{\bf c.}]
$F(\Phi_{X_k}(\tau_i))\subseteq U_i$
\item[{\bf d.}] $\Phi_{X_{k+1}}(\tau_i) = \Phi_{X_{k}}(\tau_i)\cap X_{k+1} $
if $i < k$
\item[{\bf e.}]
$\overline{U_i}\cap \overline{U_j} = \emptyset$ if $\tau_i$ and $\tau_j$ are incompatible
\item[{\bf f.}] $U_i \subseteq U_j$ if $\tau_j
<\tau_i$ \item[{\bf g.}] if $\tau_i < \tau_j$ then $\overline{U_j} \cap \overline{F(\Phi_{X_k} (\tau_i) \setminus
\Phi_{X_k}(\tau_j))} = \emptyset$
\item[{\bf h.}] $X_k$ satisfies Condition {\bf B.} for each $k \in \omega$
\end{itemize} If this can be accomplished then let $X = \cap_{k\in \omega}X_k$.  It follows that
$X$ is large and closed because, by (d),  branching is eventually preserved at each node.
Moreover $F\restriction X$ is also one-to-one
because of the choice of the $U_i$ satisfying (e) for each $i \in\omega$.  To see that $F$ is a
homeomorphism onto its range suppose that $V \subseteq X$ is an open set and that $z$ belongs to the
image of $V$ under $F$. This means that there is some $i\in \omega$ and $z'$ such that
$z' \in \Phi_X(\tau_i) \subseteq V$ and $F(z') = z$. It follows that $z \in U_i\cap F(X)$ and so it
suffices to show that $U_i\cap F(X) = F(\Phi_X(\tau_i))$.    Clearly (c) implies that
$U_i\cap F(X) \supseteq F(\Phi_X(\tau_i))$. On the other hand, if $w \in U_i\cap F(X)
$ then there
is some $w' \in X$ such that $F(w') = w$. Since $w\in U_i$ it follows  that $w' \in
\Phi_{X_k}(\tau_i)$ for each $k \geq i$ because $\{\Phi_{X_k}(\tau_j)  : j\in \omega\}$ is a tree
base.                             Hence $w \in F(\Phi_X(\tau_i))$.

To perform the  induction, use the hypothesis    on  $\{\tau_i : i\in k\}$ to choose
a maximal $\tau_i $ below $\tau_k$.
Hence there is a unique $\rho$ such that $\tau_k(\mu) = \tau_i(\mu)$ if $\mu\neq \rho$
and $\tau_k(\rho) = \tau_i(\rho)\wedge W$ for some integer $W$.
The open set $U_k$ will be chosen
so that $\overline{U_k} \subseteq U_i$ and this will guarantee that if $\tau_j$
is incompatible with $\tau_i$ then
$\overline{U_k}\cap \overline{U_j} = \emptyset$. The hypothesis on $\{\tau_i :
i\in k\}$ also implies that there is no $j \in k$ such
that $\tau_k < \tau_j$. Moreover, if $\tau_i < \tau_j$ then
$\overline{F(\Phi_{X_k}(\tau_i) \setminus \Phi_{X_k}(\tau_j))}\cap \overline{U_j} =
\emptyset$.

To satisfy
Condition (g), let $\{\delta_m : m\in a\}$   enumerate, in increasing order,
the domain of $   \tau_i$ together with the unique ordinal $\rho$ and
define $H : \presup{a}{\omega}
\to \RR$ as follows. Choose $y_s \in \presup{\alpha}{(\fomom)}$
so that for each $s \in
\presup{a}{\omega}$: \begin{itemize}
\item $y_s \in \Phi_{X_k}(\tau_i\wedge s)$ where, in this context, $\tau_i\wedge s$ is defined
by $(\tau_i \wedge s)(\delta_m) = \tau_i(\delta_m)\wedge s(m)$
\item if $s\restriction j = s'\restriction j$ then $y_s \restriction \delta_j
=
   y_{s'}\restriction \delta_j  $    \item if $s\neq s' $ then $F(y_s) \neq F(y_{s'})$
\end{itemize}
This is easily done using Condition {\bf B.} to satisfy the last two
conditions. Finally, define $H(s) = F(y_s)$ and observe that this is one-to-one.

Now use Lemma~\ref{FullSubtree} to find a full subset $T\subseteq \presup{a}{\omega}$ such that
$ H\restriction T$ has discrete image, and furthermore, this is witnessed by $\{{\cal V}_t : t\in T\}$.
Shrinking $T$ by a finite amount, if necessary, it may be assumed that
$\Phi_{X_k}(\tau_j) \cap \Phi_{X_k}(\tau_i\wedge s) = \emptyset$
 for all
$s\in T$ and $j\in k$ because $a \geq 1$.
Let $$X_{k+1} =( X_k \setminus \Phi_{X_k}(\tau_i))  \cup
(\cup\{\Phi_{X_k}(\tau_i\wedge s) : s\in T\})\cup(\cup
\{\Phi_{X_k}(\tau_j) : \tau_i \leq \tau_j\}) $$ and define
$U_k = {\cal V}_{\bar{t}}\cap U_i$ where $\bar{t}\in T$ is lexicographically the
first element of $T$. It is an easy
matter to verify that all of the induction hypotheses are satisfied.\stopproof

To finish the proof of the Lemma~\ref{keylem}   suppose that $\xi\in\omega_2+1$,
$\SuperPerfect_\xi$ is the  iteration with
countable support of the partial orders $\SuperPerfect$.
Suppose also that $p\forces{\SuperPerfect_\xi}{x\in \Cantor}$
 and
$$p\forces{\SuperPerfect_\xi}{H:\Cantor\to\Cantor  \mbox{ is a Borel function}}$$
Let $\eta\in \omega_2$ be such that $x$ occurs for
the first time in the model $V[G\cap \SuperPerfect_\eta]$.
Let ${\frak M}$ be a countable elementary submodel of $H({(2^{\aleph_0})}^+)$
containing $p$ and the names $x$ and $H$.  It follows from
Lemma~\ref{ShelCountSup} that it is
possible to find $q \leq p$ which is strongly $\PP_\eta$-generic over $\frak M$.
Let $F :  S({\frak M}\cap \xi, q) \to \Cantor$ be defined by $F(\Gamma) = x_{\Gamma}$ or,
in other words, $F(\Gamma)$ is the interpretation of $x$ in ${\frak M}[\Gamma]$.
It follows from the second clause of Definition~\ref{StronglyGeneric} that $F$ is a continuous
function. Moreover,
because it is assumed that $x$ does not belong to any model ${\frak M}[G\cap
\SuperPerfect_\mu]$ where $\mu\in\eta$, it follows that Condition
{\bf B.} of
Lemma~\ref{GenericCondition} is satisfied by $F$.  Using this lemma, and the
fact that $\eta\cap {\frak M}$ has countable order type, it is possible to find
$q' \leq q$ such that $\dom(q) = \dom (q')$ and
$F \restriction S({\frak
M}\cap\eta,q')$ is a homeomorphism onto its range.

Now let
$X$ be the image of $S({\frak M}\cap \eta,q')$ under the mapping $F$.
An inspection of the definition of $S({\frak M}\cap \eta,q')$ reveals it to be a Borel set.
Since $F\restriction S({\frak M}\cap \eta,q')$ is one-to-one, it follows that $X$ is also Borel.
Obviously
$q'\forces{\SuperPerfect_{\omega_2}}{x\in X}$.  Because the name $H$ belongs
to $\frak M$ and $F$ is one-to-one on $X$, it is possible to define a mapping
$H' :X \to \Cantor$ by defining $H'(z) $ to be the interpretation of $H(x)$ in ${\frak
M}[F^{-1}(z)]$.  Obviously $q'\forces{\SuperPerfect_{\omega_2}}{H(x) = H'(x)}$.

All that remains to be shown is that $H'$ is continuous.  To see this, let $z \in X$.
Then there is some $\Gamma \in S({\frak M}\cap \eta,q'')$   such that $z =
F(\Gamma) = x_\Gamma$. For any interval with rational end-points,  $(p,q)$,
the
statement $\psi_{p,q} $ which asserts that $H(x)\in (p,q)$
has all of its parameters in $\frak M$.  Moreover,
${\frak M}[\Gamma]\models H(x) = H(x_\Gamma) = H'(z)$.
For each interval with
rational end-points containing $H'(z)$, $(p,q)$,
there is therefore
an open neighbourood $U_{p,q}$ of $\Gamma$ such that ${\frak M}[\Gamma']\models
\psi_{p,q}$ for each $\Gamma' \in U_{p,q}$. Since $F\restriction
S({\frak M}\cap
\eta,q'')$ is a homeomorphism, it follows that the image of any $U_{p,q}$ under $F$
is an open neighbourhood $U_{p,q}^*$ of $z$. Now, if $\bar{z} \in U_{p,q}^*$ then
$\bar{z} =
x_{\Gamma'} $ for some $\Gamma'\in U_{p,q}$ and, therefore ${\frak
M}[\Gamma']\models \psi_{p,q}$. This means that the interpretation of
$H(x)$ in ${\frak M}[\Gamma']$ belongs to $(p,q)$.
Hence the image of $U_{p,q}^*$ under $H'$ is contained in $(p,q)$ and so
$H'$ is continuous.

\section{Remarks}
The proof presented here can also be generalised, without difficulty, to apply
to the iteration of $\omega_2$ Laver reals as well super-perfect reals. The
notion of
a large set has its obvious analogue which can be used to deal with the1 iteration.
In the single step case use the  proof  that a Laver real is minimal
\cite{gros.cift}.
The only difference is that, for a Laver condition $T$, the ``frontiers'' of \cite{gros.cift}
should be used in place of the images of $\theta_T\restriction
\presup{n}{\omega}$.
 In fact,
the proof of the preceding section can be viewed as a generalisation of the fact
that adding super-perfect real adds a minimal real in the sense that the
structure of the iterated model is shown to depend very predictably on the
generic reals added.

\makeatletter \renewcommand{\@biblabel}[1]{\hfill#1.}\makeatother
\renewcommand{\bysame}{\leavevmode\hbox to3em{\hrulefill}\,}

\end{document}